\documentclass[a4paper,11pt]{article}

\usepackage{amssymb,amsmath,amsthm,times,fullpage}
\usepackage{tikz}
\usepackage[pdftex]{hyperref}
\usepackage{enumitem}
\usepackage{multirow}
\usepackage{float}
\hypersetup{plainpages=True, pdfstartview=FitH, bookmarksopen=true, colorlinks=true,linkcolor=blue,citecolor=blue}
\usepackage{graphicx}

\def\pf{\noindent \emph{Proof.}\ }
\def\qed{{\quad\rule{1mm}{3mm}\,}}

\begin{document}

\setlist[description]{font=\normalfont\space}

\newtheorem{thm}{Theorem}
\newtheorem{cor}{Corollary}
\newtheorem{lmm}{Lemma}
\newtheorem{conj}{Conjecture}
\newtheorem{pro}{Proposition}
\newtheorem{Def}{Definition}
\theoremstyle{remark}\newtheorem{Rem}{Remark}

\title{Limit Theorems for Patterns in Ranked Tree-Child Networks}
\author{Michael Fuchs\thanks{Supported by MOST under the research grant MOST-109-2115-M-004-003-MY2.}, Hexuan Liu, Tsan-Cheng Yu\\
    Department of Mathematical Sciences\\
    National Chengchi University\\
    Taipei 116\\
    Taiwan}

\maketitle

\begin{abstract}
We prove limit laws for the number of occurrences of a pattern on the fringe of a ranked tree-child network which is picked uniformly at random. Our results extend the limit law for cherries proved by Bienvenu et al. (2022). For patterns of height $1$ and $2$, we show that they either occur frequently (mean is asymptotically linear and limit law is normal) or sporadically (mean is asymptotically constant and limit law is Poisson) or not all (mean tends to $0$ and limit law is degenerate). We expect that these are the only possible limit laws for any fringe pattern.
\end{abstract}

\vspace*{0.25cm}{\it Keywords:} Phylogenetic network, tree-child network, ranked tree-child network, pattern, limit law, method of moments.

\vspace*{0.25cm}{\it Mathematics subject classification (2020):} 05C20, 60C05, 60F05, 92D15.

\section{Introduction}

Studying properties of shape statistics for random models that are used to describe the evolutionary relationship between species is an important topic in biology. For {\it phylogenetic trees}, which are used to model {\it non-reticulate evolution}, many such studies have been performed and the stochastic behavior of, e.g., pattern counts are known in great detail; see \cite{ChFu,ChKaWu,ChThWu,DiWi,HoJa,KaChWu,KeSt,Ro,WuCh}. On the other hand, for {\it phylogenetic networks}, which are used to model {\it reticulate evolution}, very little is known about the number of occurrences of patterns when the networks from a given class are randomly sampled. This is due to the fact that even counting questions for phylogenetic networks from a given class were still open until recently; see \cite{CaZh,ChFuLiWaYu,DiSeWe,FuLiYu,FuYuZh1,FuYuZh2,GuRaZh,PoBa} for progress on counting questions for some of the major classes of phylogenetic networks.

One idea which made the above mentioned counting questions easier and also allowed the investigation of stochastic properties of shape parameters was the idea of ranking phylogenetic networks; see \cite{BiLaSt}. (In fact, as argued in \cite{BiLaSt}, ranked networks might be more important from a practical point of view because these networks are obtained from an evolution process; see the definition below.) More precisely, the author of \cite{BiLaSt} defined {\it ranked tree-child networks} and proved several results for them, e.g., they obtained a Poisson limit law for the number of cherries when a network is sampled uniformly at random. This is, as far as we know, the first limit result for the number of occurrences of a pattern in phylogenetic networks. (Note that cherries were also the first pattern for which a limit law was proved for phylogenetic trees; see \cite{KeSt}.)

The main goal of this work is to study limit laws for the number of occurrences of other patterns on the fringe of ranked tree-child networks. Our results might give an indication of what to expect when studying patterns for other classes of phylogenetic networks.

Before explaining our results, we will give definitions and recall previous results. We start with the definition of a (rooted, binary) phylogenetic network which is a directed acyclic graph (DAG) without double edges such that every node falls into one of the following four categories:
\begin{itemize}
\item A (unique) {\it root} which has in-degree $0$ and out-degree $1$;
\item {\it Leaves} which have in-degree $1$ and out-degree $0$ and which are bijectively labeled by $\{1,\ldots,n\}$ where $n$ is their number;
\item {\it Tree nodes} which are nodes of in-degree $1$ and out-degree $2$;
\item {\it Reticulation nodes} which are nodes of in-degree $2$ and out-degree $1$.
\end{itemize}

\begin{figure}[t]
    \centering
    \includegraphics[scale=0.9]{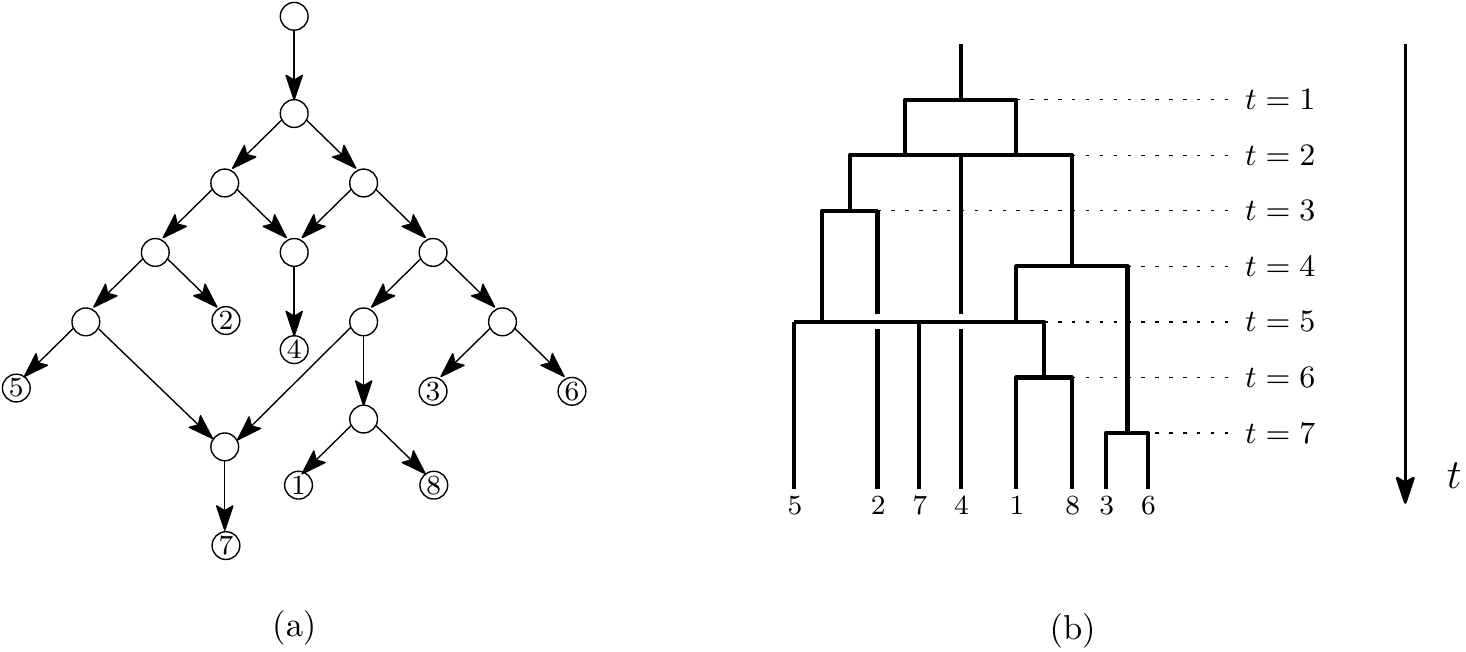}
    \caption{(a) A phylogenetic network which satisfies the tree-child property and is rankable; (b) The rankable tree-child network from (a) with a ranking.}\label{RTCN+ranking}
\end{figure}

See Figure~\ref{RTCN+ranking}-(a) for an example. In fact, the network from this figure is even a {\it tree-child network}.

\begin{Def}
A phylogenetic network is called tree-child network if every non-leaf node has at least one child which is not a reticulation node.
\end{Def}

\begin{figure}[h!]
\centering
\includegraphics[scale=0.9]{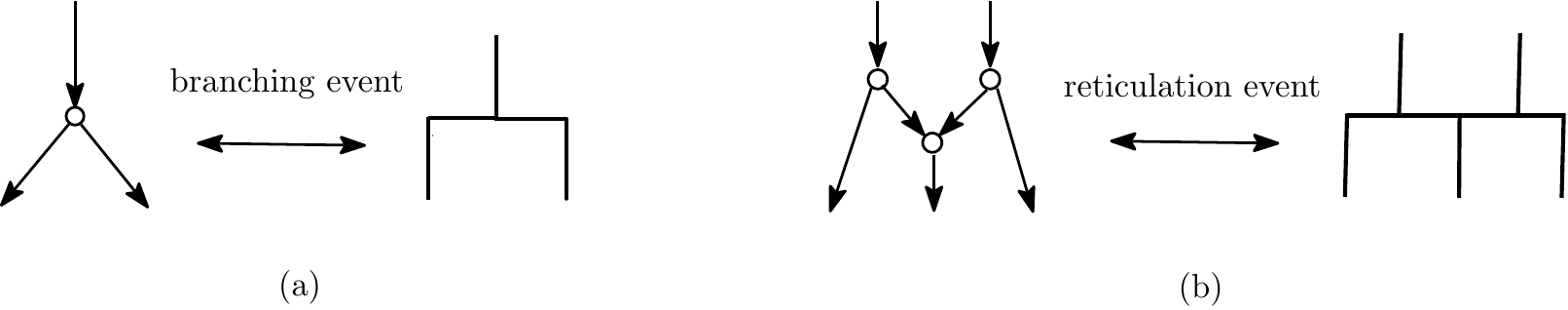}
\caption{The branching and reticulation event used in the construction of ranked tree-child networks.}\label{bra-ret}
\end{figure}

For a tree-child network, we call a tree-node a {\it branching event} and a reticulation node with its two parents a {\it reticulation event}; see Figure~\ref{bra-ret} for the graphical depiction of these two events that we are going to use in the sequel. (Vertical edges in this depiction will subsequently be called {\it lineages}.)

\vspace*{0.35cm}

The main object in this paper are ranked tree-child networks which will be defined next.

\begin{Def}
A tree-child network is called rankable if it has recursively evolved starting from a branching event by attaching in each step either a branching event or a reticulation event. A rankable tree-child network together with a ranking of its events is called a ranked tree-child network.
\end{Def}

See Figure~\ref{RTCN+ranking}-(b) for a ranked tree-child arising from the rankable tree-child network from Figure~\ref{RTCN+ranking}-(a).

Ranked tree-child networks have been introduced in \cite{BiLaSt} where the authors proved a wealth of combinatorial and stochastic properties about them. Two of these properties concern the number of occurrences of patterns in a ranked tree-child network which is picked uniformly at random from all ranked tree-child networks with $n$ leaves. (We call such a network a {\it random ranked tree-child network} in the sequel.) To recall these results, we need two definitions. First, a {\it cherry} is a tree node with both children leaves (or equivalently, a branching event with both outgoing lineages external); a {\it trident} is a branching event with all three outgoing lineages external. For instance, in the ranked tree-child network from Figure~\ref{RTCN+ranking}-(b), there are $2$ cherries and no trident.

Denote by $C_n$ resp. $T_n$ the number of cherries resp. tridents in a random ranked tree-child network with $n$ leaves. The following limit law result was proved in \cite{BiLaSt} for $C_n$: as $n\rightarrow\infty$, $C_n$ weakly converges to the Poisson distribution with parameter $1/4$, i.e.,
\begin{equation}\label{ll-cherries}
C_n\stackrel{d}{\longrightarrow}{\rm Poisson}(1/4),\qquad (n\rightarrow\infty).
\end{equation}
On the other hand, for the number of tridents, the authors in \cite{BiLaSt} just proved a weak law of large numbers:
\begin{equation}\label{wlln}
\frac{T_n}{n}\stackrel{{\mathbb P}}{\longrightarrow}\frac{1}{7},\qquad (n\rightarrow\infty).
\end{equation}
Our first result improves this to a central limit theorem.
\begin{thm}\label{main-result-1}
For the number of tridents $T_n$ in a random ranked tree-child network with $n$ leaves, we have
\begin{equation}\label{ll-tridents}
\frac{T_n-n/7}{\sqrt{24n/637}}\stackrel{d}{\longrightarrow} N(0,1),\qquad (n\rightarrow\infty),
\end{equation}
where $N(0,1)$ denotes the standard normal distribution.
\end{thm}

\begin{figure}[t]
    \centering
    \includegraphics[scale = 0.9]{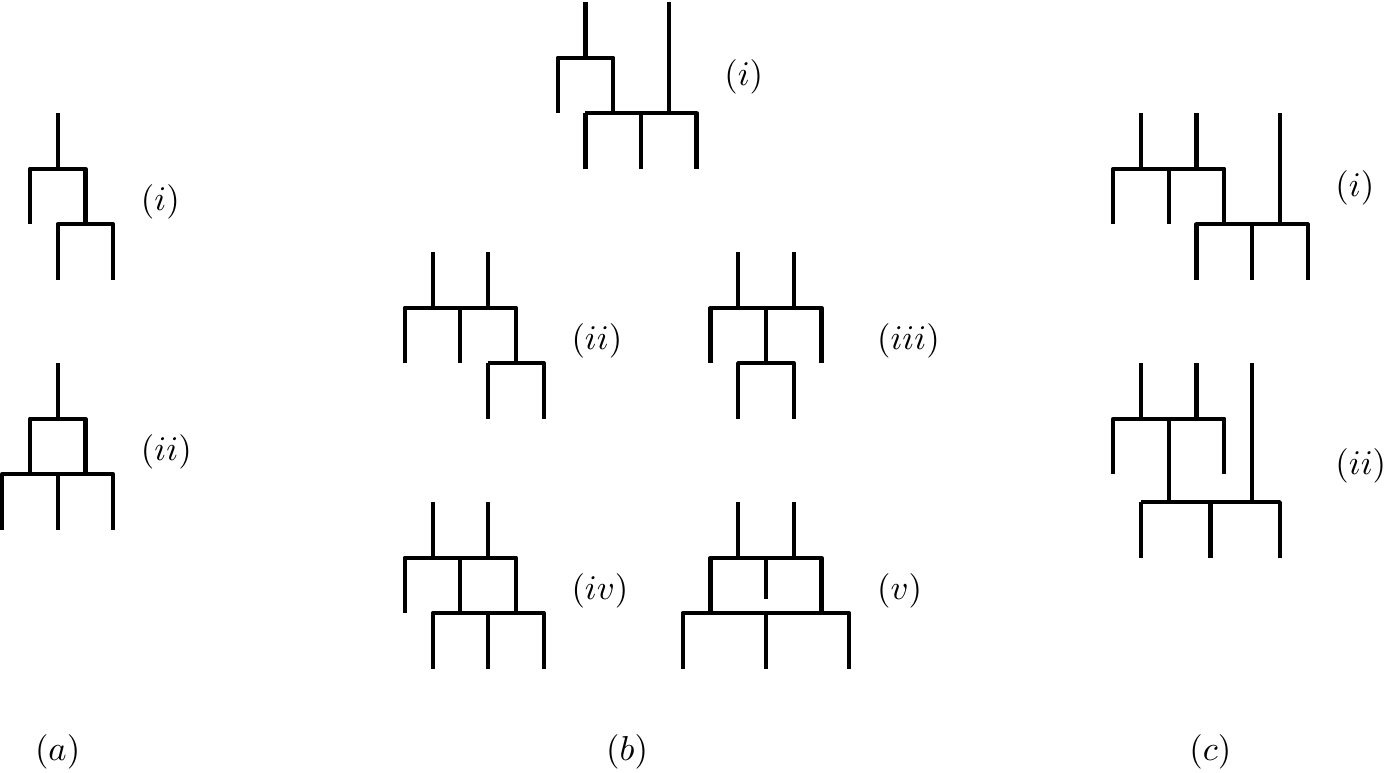}
    \caption{All patterns of height 2. The number of occurrences of these patterns in a ranked tree-child networks with a large number of leaves is as follows: (a) These two do not occur; (b) These five occur only sporadically; (c) These two occur frequently. (For details see Theorem~\ref{main-result-2}.)}\label{pattern-2}
\end{figure}

Note that (\ref{ll-cherries}) and (\ref{ll-tridents}) are all the limit laws of  patterns of height $1$, where the height is defined as the number of steps in the evolution process from the definition of ranked tree-child networks. We next turn to all patterns of height $2$ which are listed in Figure~\ref{pattern-2}. For the number of occurrences of these patterns in random ranked tree-child networks with $n$ leaves, we have the following result.

\begin{thm}\label{main-result-2}
Denote by $X_n$ the number of occurrences of a (fixed) pattern of height $2$  in a random ranked tree-child networks with $n$ leaves. Then, we have the following limit law results.
\begin{description}
\item[(A)] For the patterns in Figure~\ref{pattern-2}-(a), we have that the limit law of $X_n$ is degenerate. More precisely,
\[
X_n\stackrel{L_1}{\longrightarrow}0,\qquad (n\rightarrow\infty).
\]
\item[(B)] For the patterns in Figure~\ref{pattern-2}-(b), we have
\[
X_n\stackrel{d}{\longrightarrow}{\rm Poisson}(\lambda),\qquad (n\rightarrow\infty),
\]
where

\vspace*{-0.1cm}
\begin{center}
\begin{tabular}{c|c|c|c|c|c}
   & (b-i) & (b-ii) & (b-iii) & (b-iv) & (b-v) \\ \hline\hline
   $\lambda$ & $1/8$ & $1/28$ & $1/56$ & $1/14$ & $1/28$
\end{tabular}
\end{center}

\item[(C)] For the patterns in Figure~\ref{pattern-2}-(c), we have
\[
\frac{X_n-\mu n}{\sigma\sqrt{n}}\stackrel{d}{\longrightarrow}N(0,1),\qquad (n\rightarrow\infty),
\]
where $(\mu,\sigma^2)=(4/77,4575916/137582445)$ and $(\mu,\sigma^2)=(2/77,2930764/137582445)$ for the patterns from Figure~\ref{pattern-2}-(c-i) and Figure~\ref{pattern-2}-(c-ii), respectively.
\end{description}
\end{thm}

Thus, the only possible limit laws for patterns of height $1$ and height $2$ are either the normal law or the Poisson law or the degenerate law. In fact, we believe that these are the only possible limit laws for the occurrence of any (fixed) pattern; see the discussion in Section~\ref{con}. Note that is a very different behavior from the one observed for patterns in random phylogenetic trees which are all known to be asymptotically normal distributed; see \cite{ChFu}.

We conclude the introduction by (briefly) explaining the proof of the above two theorems and at the same time presenting an outline of the paper. The proofs use the same strategy that was used for (\ref{ll-cherries}) in \cite{BiLaSt}, namely, coupling with a Markov chain and applying the method of moments. In addition, an important step in the proof of Theorem~\ref{main-result-1} (as well as in the central limit theorems of Theorem~\ref{main-result-2}) will be {\it shifting-the-mean} which will considerably reduce the complexity of computing higher central moments. The proof of Theorem~\ref{main-result-1} will be presented in the next section. The proof of Theorem~\ref{main-result-2}, which will be given in Section~\ref{height-2}, will be more demanding since for all patterns from Figure~\ref{pattern-2} except for the cases (b-i), (c-i) and (c-ii), we will need to consider the pattern together with another pattern to be able to compute moments; for the remaining cases, namely, Figure~\ref{pattern-2}-(b-i), (c-i) and (c-ii), we will need to consider each of them with two more patterns (one of which will be a pattern of height $3$) and the description of the Markov chain will require many cases; this will make the computation of the moments complicated (and thus the computations will be done with the help of Maple). The paper will be concluded with some remarks in Section~\ref{con}.

\section{Patterns of Height 1}

In this section, we are going to prove Theorem~\ref{main-result-1}. Before doing so, we will need some preliminaries.

First, we recall a random process from \cite{BiLaSt} (which was called {\it forward construction} in that paper). Start with a branching event and recursively attach in each consecutive step either a branching or a reticulation event. This is done as follows: in the $\ell-1$-st step a tuple $(\ell_1,\ell_2)$ of the $\ell$ external lineages is picked uniformly at random and a branching event is attached to the lineage $\ell_1$ if $\ell_1=\ell_2$ or a reticulation event with incoming lineages $\ell_1$ and $\ell_2$ is created if $\ell_1\ne \ell_2$. The process stops if $n$ external lineages are created. Note that the resulting network is {\it not} a ranked tree-child network since (a) leaves are not labeled and (b) an order of the incoming lineages of a reticulation event has been fixed. However, it was proved in \cite{BiLaSt} that the random variable counting the number of tridents (or more generally any pattern) in the network resulting from the the above process has the same distribution as the corresponding random variable in random ranked tree-child networks.

By analyzing what will happen with the number of tridents by performing one more step in the above construction, the authors in \cite{BiLaSt} obtained the following result.

\begin{lmm}[Bienvenu et al.; \cite{BiLaSt}]
Define a Markov process by $T_2=0$ and
\begin{equation}\label{markov-Tn}
(T_{n+1}\vert T_n=j)=\begin{cases} j-1,&\text{with probability}\ {\displaystyle \frac{3j(3j-2)}{n^2}};\vspace*{0.15cm}\\
j,&\text{with probability}\  {1-\displaystyle\frac{3j(3j-2)+(n-3j)(n-3j-1)}{n^2}};\vspace*{0.15cm}\\j+1,&\text{with probability}\ {\displaystyle\frac{(n-3j)(n-3j-1)}{n^2}}.
\end{cases}
\end{equation}
Then, $T_n$ has the same distribution as the number of tridents in a random ranked tree-child network with $n$ leaves.
\end{lmm}

This was then used to compute the mean and variance of $T_n$ (which then in turn implied (\ref{wlln})). Since we will need the result for the mean, we recall it below.

\begin{lmm}[Bienvenu et al.; \cite{BiLaSt}]\label{mean-BiLaSt}
Let $\mu_n:={\mathbb E}(T_n)$. Then, $m_n$ satisfies the recurrence
\[
\mu_{n+1}=\left(1-\frac{3}{n}\right)^2\mu_n+1-\frac{1}{n}
\]
whose solution is given by
\begin{equation}\label{exp-mu}
\mu_n=\frac{(15n^3-85n^2+144n-71)n}{105(n-1)(n-2)(n-3)},\qquad (n\geq 4).
\end{equation}
\end{lmm}

Note that this result gives that $\mu_n\sim n/7$ as $n\rightarrow\infty$. An asymptotic expansion for the variance was also obtained in \cite{BiLaSt} by first computing the second moment (which also satisfies a similar first-order recurrence as the mean) and then using the definition of the variance. It is actually easier to work directly with the variance which also satisfies a similar recurrence as the mean. (The reason why this is easier is that it already incorporates the cancellations which arise from ${\mathbb E}(T_n^2)\sim({\mathbb E}(T_n))^2$; this method is the {\it shifting-the-mean} method mentioned in the last paragraph of the previous section.) Moreover, this also extends to higher moments.

More precisely, we define
\[
\phi_{n,m}:={\mathbb E}(T_n-\mu_n)^{m}.
\]
Then, this sequence satisfies the following recurrence.
\begin{lmm}\label{rec-cent-mom}
For $m\geq 2$, we have
\[
\phi_{n+1,m}=\left(1-\frac{3m}{n}\right)^2\phi_{n,m}+\psi_{n,m}
\]
with
\[
\psi_{n,m}=\sum_{j=0}^{m-1}\phi_{n,j}\Lambda_j(n),
\]
where $\Lambda_j(n)$ admits the complete asymptotic expansion
\[
\Lambda_j(n)\sim\sum_{\ell=0}^{\infty}\frac{\lambda_{j,\ell}}{n^\ell},\qquad (n\rightarrow\infty)
\]
with $\lambda_{j,\ell}\in{\mathbb R}$; in particular, $\lambda_{m-1,0}=0$ and
\[
\lambda_{m-2,0}=\binom{m}{2}\times\frac{24}{49}.
\]
\end{lmm}

\pf Set $\bar{T}_n:=T_n-\mu_n$. From (\ref{markov-Tn}), we have
\begin{align*}
\phi_{n+1,m}
=&\ {\mathbb E}(\bar{T}_{n}+\mu_n-\mu_{n+1}-1)^m\frac{3(\bar{T}_n+\mu_n)(3(\bar{T}_n+\mu_n)-2)}{n^2}\\
&+{\mathbb E}(\bar{T}_{n}+\mu_n-\mu_{n+1})^m\\
&\qquad\times\left(1-\frac{3(\bar{T}_n+\mu_n)(3(\bar{T}_n+\mu_n)-2)+(n-3(\bar{T}_n+\mu_n))(n-3(\bar{T}_n+\mu_n)-1)}{n^2}\right)\\
&+{\mathbb E}(\bar{T}_{n}+\mu_n-\mu_{n+1}+1)^m\frac{(n-3(\bar{T}_n+\mu_n))(n-3(\bar{T}_n+\mu_n)-1)}{n^2}.
\end{align*}
From this, by expanding what is inside the means by the binomial theorem, we see that
\[
\psi_{n,m}=\sum_{j=0}^{m+2}\phi_{n,j}\Lambda_j(n).
\]
However, $\Lambda_{m+2}(n)=0$ (because the probabilities in (\ref{markov-Tn}) sum up to $1$) and straightforward computation (best done with a computer algebra system such as Maple) shows that $\Lambda_{m+1}(n)=\Lambda_{m}(n)=0$. Next, since Lemma~\ref{mean-BiLaSt} implies that $\mu_n-\mu_{n+1}={\mathcal O}(1)$, we have
\[
\Lambda_j(n)\sim\sum_{\ell=0}^{\infty}\frac{\lambda_{j,\ell}}{n^\ell},\qquad (n\rightarrow\infty)
\]
and by some more computations (again best done with Maple), we obtain that $\lambda_{m-1,0}=0$ and that $\lambda_{m-2,0}$ is as claimed.\qed

Note that the above recurrence has the (general) form
\begin{equation}\label{gen-rec}
\phi_{n+1}=\left(1-\frac{\kappa}{n}\right)^2\phi_n+\psi_n,\qquad (n\geq\kappa+1)
\end{equation}
with a suitable initial value $\phi_{\kappa+1}$, where $\kappa\in{\mathbb N}$ and $\{\psi_n\}_{n\geq \kappa+1}$ is a given sequence. We need a general result for such a sequence. (This result was also implicitly contained in \cite{BiLaSt}.)

\begin{lmm}\label{asymp-transfer}
Assume that $\phi_n$ satisfies (\ref{gen-rec}). If $\psi_n\sim cn^{\alpha}$ with $\alpha>-2\kappa-1$ a real number, then
\[
\phi_n\sim\frac{c}{2\kappa+\alpha+1}n^{\alpha+1},\qquad (n\rightarrow\infty).
\]
\end{lmm}
\begin{Rem}
Throughout the paper, the notation $f(n)\sim cg(n)$ (with $g(n)>0$) means that $f(n)=cg(n)+o(g(n))$. (Note that this is the usual meaning if $c\ne 0$; however, if $c=0$, then $f(n)=o(g(n)$.)
\end{Rem}

\pf Iterating (\ref{gen-rec}) gives the solution
\begin{equation}\label{sol-gen-rec}
\phi_n=\binom{n-1}{\kappa}^{-2}\left(\psi_{\kappa+1}+\sum_{\ell=\kappa+1}^{n-1}\binom{\ell}{\kappa}^2\psi_{\ell}\right)
\end{equation}
Note that
\begin{equation}\label{binom-asym}
\binom{\ell}{\kappa}^2\sim\frac{\ell^{2\kappa}}{\kappa!^2},\qquad (\ell\rightarrow\infty).
\end{equation}
Thus, if $\psi_n\sim cn^{\alpha}$, then
\[
\psi_{\kappa+1}+\sum_{\ell=\kappa+1}^{n-1}\binom{\ell}{\kappa}^2\psi_{\ell}\sim\frac{c}{\kappa!^2}\sum_{\ell=\kappa+1}^{n-1}\ell^{2\kappa+\alpha}\sim\frac{c}{\kappa!^2(2\kappa+\alpha+1)}n^{2\kappa+\alpha+1}.
\]
Plugging this into (\ref{sol-gen-rec}) and using once more (\ref{binom-asym}) gives the claimed result.\qed

Applying the last result to the recurrence for the central moments (Lemma~\ref{rec-cent-mom}) and using induction gives the following asymptotic result for all central moments of $T_n$. (This method is sometimes refered to as {\it moment-pumping} in the literature; see Section VII.10.1 in \cite{FlSe}.)

\begin{pro}\label{mom-Tn}
As $n\rightarrow\infty$, the central moments of $T_n$ satisfy
\[
{\mathbb E}\left((T_n-\mu_n)^m\right)\sim g_{m}\left(\frac{24}{637}\right)^{m/2}n^{m/2}.
\]
Here, $g_m$ denotes the $m$-th moment of the standard normal distribution, i.e.,
\begin{equation}\label{moment-N01}
g_m={\mathbb E}\left(N(0,1)^m\right)=\begin{cases}{\displaystyle\frac{m!}{2^{m/2}(m/2)!}},&\text{if}\ m\ \text{is even};\\0,&\text{if}\ m\ \text{is odd}.\end{cases}
\end{equation}
\end{pro}

\pf The claim (trivially) holds for $m=0$ and $m=1$. Next, we assume that the claim holds for all $m'<m$. In order to show it for $m$, note that from the induction hypothesis, we have for $\psi_{n,m}$ from Lemma~\ref{rec-cent-mom}:
\[
\psi_{n,m}\sim\binom{m}{2}\times\frac{24}{49}\times g_{m-2}\left(\frac{24}{637}\right)^{m/2-1}n^{m/2-1}.
\]
Applying Lemma~\ref{asymp-transfer}, we obtain that
\[
\phi_{n,m}\sim\binom{m}{2}\times\frac{24}{49}\times \frac{g_{m-2}}{6m+m/2}\left(\frac{24}{637}\right)^{m/2-1}n^{m/2}\sim (m-1)g_{m-2}\left(\frac{24}{637}\right)^{m/2}n^{m/2}
\]
from which the claimed result follows since $(m-1)g_{m-2}=g_{m}$.\qed

Theorem~\ref{main-result-1} follows now from the last proposition by the Fr\'{e}chet-Shohat Theorem; see, e.g., Section~5.30 in \cite{Bi}.

\section{Patterns of Height 2}\label{height-2}

In this section, we will prove Theorem~\ref{main-result-2}. This will be done in the three paragraphs below, one for each of the three cases in Theorem~\ref{main-result-2}. The proof will proceed along similar lines as the proof from the previous section. In addition, it will make use of the expansions from Proposition~\ref{mom-Tn} which imply that
\begin{equation}\label{exp-mom-Tn}
{\mathbb E}(T_n^{m})\sim\mu_n^{m}\sim\frac{n^{m}}{7^m},\qquad (n\rightarrow\infty)
\end{equation}
and
\begin{equation}\label{exp-mom-Cn}
{\mathbb E}(C_n^{\underline{m}}):={\mathbb E}(C_n(C_n-1)\cdots(C_n-m+1))\sim\frac{1}{4^m},\qquad (n\rightarrow\infty)
\end{equation}
which was used in \cite{BiLaSt} to prove (\ref{ll-cherries}).

\paragraph{Degenerate Limit Laws.} We will only consider the pattern in Figure~\ref{height-2}-(a-i); the other pattern whose limit law is degenerate, namely the one in Figure~\ref{height-2}-(a-ii), is treated similarly. Our method below will be slightly more general than needed; the reason for this is that the other cases from Theorem~\ref{main-result-2} will be proved with similar arguments (and in these cases, the generality below is needed).

We will consider the pattern in Figure~\ref{height-2}-(a-i) (called a pattern of type $A$ in the sequel) together with a cherry which is not contained in a pattern of type $A$ (called a pattern of type $B$ in the sequel). Assume that a random ranked tree-child network with $n$ leaves contains $a$ pattern of type $A$ and $b$ pattern of type $B$, respectively. Note that each external lineage belongs either exactly to a pattern of type $A$ or a pattern of type $B$ or to neither of these patterns (such a lineage will be called a pattern of type $C$ in the sequel). We now carefully list what will happen with the number of patterns of type $A$ and type $B$ if either a branching or a reticulation event is added in the next step of the forward construction that was described at the beginning of the last section.

First, if a branching event is added, then we have the following cases.

\vspace{0.1cm}
\begin{table}[!h]
\begin{center}
\begin{tabular}{c|c|c|c}
   & type $A$ & type $B$ & probability \\ \hline\hline
   \multirow{2}{*}{type $A$} & $-1$ & $+2$ & $a/n^2$ \\
   & $0$ & $0$ & $2a/n^2$ \\ \hline
   type $B$ & $+1$ & $-1$ & $2b/n^2$ \\ \hline
   type $C$ & $0$ & $+1$ & $(n-3a-2b)/n^2$
\end{tabular}
\end{center}
\vspace{-0.2cm}
\caption{The change of the number of patterns of type $A$ and type $B$ (first and second column) if a branching event is attached to an external lineage belonging to a pattern of type $A$, type $B$ or type $C$ (rows). Here, $a$ and $b$ denote the number of patterns of type $A$ and $B$, respectively. The probability for each cases is given in the third column if one starts with a ranked tree-child network with $n$ leaves (and thus the probability of picking a fixed external lineage is $1/n^2$).}\label{bran-(a-i)}
\end{table}

For instance, if the branching event is attached to an external lineage from a pattern of type $A$ and that lineage does not belong to the cherry in that pattern, then two patterns of type $B$ are created whereas one pattern of type $A$ was destroyed; this is the first sub-row of the first row in the above table and the probability that this happens is given by the number of possible choices of such a lineage ($a$) divided by the number of choices of a pair of two external lineages ($n^2$). Likewise, the second sub-row of the first row is the case where the branching event is attached to an external lineage from a pattern of type $A$ where this lineage now belongs to the cherry of that pattern. Similarly, the remaining rows are explained.

Next, if a reticulation event is added, then we have the cases listed in Table~\ref{ret-(a-i)}.

\begin{table}[t]
\begin{center}
\begin{tabular}{c|c|c|c}
   & type $A$ & type $B$ & probability \\ \hline\hline
   $A$ & $-1$ & $0$ & $6a/n^2$ \\ \hline
   \multirow{3}{*}{$A\ \&\ A$} & $-2$ & $0$ & $4a(a-1)/n^2$\\
   & $-2$ & $+1$ & $4a(a-1)/n^2$ \\
   & $-2$ & $+2$ & $a(a-1)/n^2$ \\ \hline
   $B$ & $0$ & $-1$ & $2b/n^2$ \\ \hline
   $B\ \&\ B$ & $0$ & $-2$ & $4b(b-1)/n^2$\\ \hline
   $C\ \&\ C$ & $0$ & $0$ & $(n-3a-2b)(n-3a-2b-1)/n^2$\\ \hline
   \multirow{2}{*}{$A\ \&\ B$} & $-1$ & $-1$ & $8ab/n^2$\\
   & $-1$ & $0$ & $4ab/n^2$ \\ \hline
   \multirow{2}{*}{$A\ \&\ C$} & $-1$ & $0$ & $4a(n-3a-2b)/n^2$\\
   & $-1$ & $+1$ & $2a(n-3a-2b)/n^2$\\ \hline
   $B\ \&\ C$ & $0$ & $-1$ & $4b(n-3a-2b)/n^2$
\end{tabular}
\end{center}
\vspace*{-.2cm}
\caption{Columns: The change of the number of patterns of type $A$ (first column) and type $B$ (second column); the probabilities of these changes are contained in the last column. Rows: a single letter means that both external lineages are picked from a pattern of that type; two letters mean that the external lineages are chosen from (different) patterns of these two types.}\label{ret-(a-i)}
\end{table}

For instance, the first sub-row of the sixth row of this column is explained as follows: if the reticulation event is attached to an external lineage belonging to the cherry of a pattern of type $A$ and an external lineage from a pattern of type $B$, then a pattern of type $A$ and type $B$ is destroyed (and no new pattern of type $B$ is created as will happen if the third external lineage from the pattern of type $A$ is chosen). The probability that this will happen is given by the number of choices of the external lineages divided by $n^2$; the number of choices equals $2\cdot 2\cdot 2\cdot ab$ since once the patterns of type $A$ and $B$ are chosen ($ab$ choices) there are $2$ choices for the external lineages in the pattern of type $A$ and $2$ in the pattern of type $B$; moreover, the factor $2$ comes from symmetry. Similarly, the other rows are explained.

Note that Table~\ref{bran-(a-i)} and Table~\ref{ret-(a-i)} give the transition probabilities of the Markov chain $(X_n,\tilde{C}_n)$ where $X_n$ and $\tilde{C}_n$ are the number of patterns of type $A$ and type $B$, respectively, in a random ranked tree-child networks of $n$ leaves. With these probabilities, we obtain the following result.

\begin{lmm}
Let $X_n$ be the number of occurrences of the pattern from Figure~\ref{height-2}-(a-i) in a random ranked tree-child network with $n$ leaves. Then,
\begin{equation}\label{rec-(a-i)}
{\mathbb E}(X_{n+1})=\left(1-\frac{3}{n}\right)^2{\mathbb E}(X_n)+\frac{2{\mathbb E}(C_n)}{n^2}.
\end{equation}
\end{lmm}
\pf For $X_n$ and $\tilde{C}_n$ (see the paragraph proceeding the lemma), we have
\begin{equation}\label{rec-Xn-Yn}
{\mathbb E}(X_{n+1})=\left(1-\frac{6}{n}+\frac{11}{n^2}\right){\mathbb E}(X_n)+\frac{2{\mathbb E}(\tilde{C}_n)}{n^2}.
\end{equation}
Note that $C_n=X_n+\tilde{C}_n$. Consequently, we obtain the claimed recurrence from the above one by replacing $\tilde{C}_n=C_n-X_n$.\qed

\begin{Rem}
Note that for (\ref{rec-Xn-Yn}), we in fact only need the first columns of Table~\ref{bran-(a-i)} and Table~\ref{ret-(a-i)}. However, all columns will be needed for the patterns below.
\end{Rem}

The last lemma now implies the limit law result from Theorem~\ref{main-result-2}-(a).

\vspace*{0.35cm}
\noindent{\it Proof of Theorem~\ref{main-result-2}-(a) for the pattern in Figure~\ref{height-2}-(a-i).}  Recurrence (\ref{rec-(a-i)}) has the form (\ref{gen-rec}) with
\[
\psi_n=\frac{2{\mathbb E}(C_n)}{n^2}\sim\frac{1}{2n^2},
\]
where the last step follows from (\ref{exp-mom-Cn}). Applying Lemma~\ref{asymp-transfer} gives
\[
{\mathbb E}(X_n)\sim\frac{1}{10n}
\]
which implies the claimed result.\qed

\paragraph{Poisson Limit Laws.} Here, we will consider the patterns from Figure~\ref{height-2}-(b). Since the proof for all patterns except for the pattern from Figure~\ref{height-2}-(b-i) is the same, we will only give details for one of these four patterns, namely, the pattern in Figure~\ref{height-2}-(b-iv).

\begin{figure}[t]
    \centering
    \includegraphics[scale=0.9]{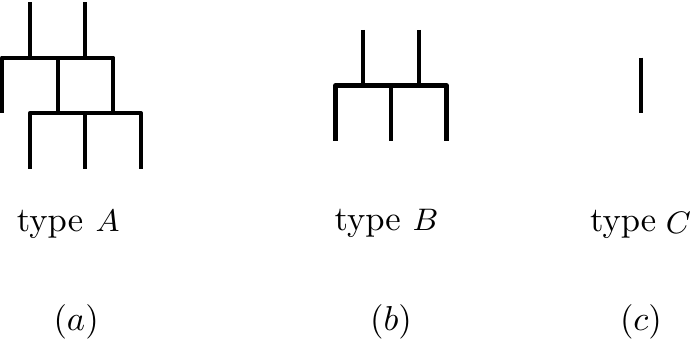}
    \caption{(a) The pattern from Figure~\ref{height-2}-(b-iv); (b) A trident which is not contained in a pattern from (a); (c) The remaining external lineages.}\label{pattern-(b-iv)}
\end{figure}

Similarly to the treatment of the pattern from the previous paragraph, we will use the types of patterns depicted in Figure~\ref{pattern-(b-iv)}. Note that each external lineage from a ranked tree-child network belongs to exactly one of these types.

Assume now that a ranked tree-child network contains $a$ patterns of type $A$ and $b$ patterns of type $B$. Then, we again consider the changes and corresponding probabilities of these numbers when a branching event or reticulation event is added in the forward construction; see Table~\ref{bra-(b-iv)} and Table~\ref{ret-(b-iv)}. This again gives the transition probabilities of a Markov chain which can be used to compute mixed moments.

\vspace*{0.1cm}
\begin{table}[!h]
\begin{center}
\begin{tabular}{c|c|c|c}
   & type $A$ & type $B$ & probability \\ \hline\hline
   \multirow{2}{*}{$A$} & $-1$ & $0$ & $3a/n^2$ \\
   & $-1$  & $+1$ & $a/n^2$ \\\hline
   $B$ & $0$ & $-1$ & $3b/n^2$ \\ \hline
   $C$ & $0$ & $0$ & $(n-4a-3b)/n^2$
\end{tabular}
\end{center}
\vspace*{-.2cm}
\caption{The change of the number of patterns of type $A$ and type $B$ (see Figure~\ref{pattern-(b-iv)}) when the next event is a branching event.}\label{bra-(b-iv)}
\end{table}

\begin{table}[!h]
\begin{center}
\begin{tabular}{c|c|c|c}
   & type $A$ & type $B$ & probability \\ \hline\hline
   \multirow{2}{*}{$A$} & $-1$ & $+1$ & $8a/n^2$ \\
   & $0$  & $0$ & $4a/n^2$ \\\hline
   \multirow{3}{*}{$A\ \&\ A$} & $-2$ & $+1$ & $9a(a-1)/n^2$\\
   & $-2$ & $+2$ & $6a(a-1)/n^2$ \\
   & $-2$ & $+3$ & $a(a-1)/n^2$ \\ \hline
   \multirow{2}{*}{$B$} & $0$ & $0$ & $2b/n^2$ \\
   & $+1$ & $-1$ & $4b/n^2$\\ \hline
   $B\ \&\ B$ & $0$ & $-1$ & $9b(b-1)/n^2$\\ \hline
   $C\ \&\ C$ & $0$ & $+1$ & $(n-4a-3b)(n-4a-3b-1)/n^2$\\ \hline
   \multirow{2}{*}{$A\ \&\ B$} & $-1$ & $0$ & $18ab/n^2$\\
   & $-1$ & $+1$ & $6ab/n^2$ \\ \hline
   \multirow{2}{*}{$A\ \&\ C$} & $-1$ & $+1$ & $6a(n-4a-3b)/n^2$\\
   & $-1$ & $+2$ & $2a(n-4a-3b)/n^2$\\ \hline
   $B\ \&\ C$ & $0$ & $0$ & $6b(n-4a-3b)/n^2$
\end{tabular}
\end{center}
\vspace*{-.2cm}
\caption{The change of the number of patterns of type $A$ and type $B$ (see Figure~\ref{pattern-(b-iv)}) when the next event is a reticulation event.}\label{ret-(b-iv)}
\end{table}

\begin{lmm} Let $X_n$ be the number of occurrences of the pattern from Figure~\ref{height-2}-(b-iv) in a random ranked tree-child network with $n$ leaves. Then, for all $r,s\geq 0$, we have
\[
{\mathbb E}(X_{n+1}^{\underline{r}}T_{n+1}^s)=\left(1-\frac{4r+3s}{n}\right)^{2}{\mathbb E}(X_n^{\underline{r}}T_n^{s})+\frac{4r}{n^2}{\mathbb E}(X_n^{\underline{r-1}} T_n^{s+1})-\frac{4r(r-1)}{n^2}{\mathbb E}(X_n^{\underline{r-1}}T_n^s)+R_n,
\]
where
\begin{equation*}
R_n=\sum_{j=0}^{s-1}\left(\lambda_{j,0}+\frac{\lambda_{j,1}}{n}+\frac{\lambda_{j,2}}{n^2}\right){\mathbb E}(X_n^{\underline{r}}T_n^{j}).
\end{equation*}
Here, $\lambda_{j,i}$ are constants which only depend on $r$ and $s$ with $\lambda_{s-1,0}=s$.
\end{lmm}

\pf Assume that $X_n=c$ and $T_n=d$. Then, we have $c$ patterns of type $A$ and $d-c$ patterns of type $B$. Using the cases from Table~\ref{bra-(b-iv)} and Table~\ref{ret-(b-iv)}, we can write $X_{n+1}^{\underline{r}}T_{n+1}^s$ given $X_n=c$ and $T_n=d$ as a sum of terms $(X_n+k)^{\underline{r}}(T_{n}+\ell)^s$ which are multiplied with the probabilities of every case and where $k\in\{-2,-1,0,1\}$ is the value from the first column and $\ell\in\{-1,0,1\}$ is the sum of the values of the first and second column for every case. Next, we replace $(X_n+k)^{\underline{r}}$ for $k=-2$ and $k=-1$ by
\[
(X_n-2)^{\underline{r}}=\frac{(X_n-r)(X_n-r-1)}{X_n(X_n-1)}X_n^{\underline{r}},\qquad(X_n-1)^{\underline{r}}=\frac{X_n-r}{X_n}X_n^{\underline{r}}
\]
and for $k=1$ by
\[
(X_n+1)^{\underline{r}}=\frac{X_n+1}{X_n-r+1}X_n^{\underline{r}}.
\]
Finally, we expand $(T_n+\ell)^s$ by the binomial theorem and simplify the sum of all terms of the same order in this expansion with Maple. (Note that since $\ell\in\{-1,0,1\}$ this only has to be done for the first three terms in the expansion; the other terms are the same only multiplied with different weights). This proves the claimed result.\qed

We are now ready to prove the Poisson limit result for the pattern from Figure~\ref{height-2}-(b-iv).

\vspace*{0.35cm}
\noindent{\it Proof of Theorem~\ref{main-result-2}-(b) for the pattern in Figure~\ref{height-2}-(b-iv).} As above, let $X_n$ denote the number of occurrences of the pattern from Figure~\ref{height-2}-(b-iv) in a random ranked tree-child network with $n$ leaves.  We will use induction to show that for all $r,s\geq 0$:
\begin{equation}\label{main-claim}
{\mathbb E}(X_n^{\underline{r}}T_n^s)\sim\frac{n^s}{14^r7^s},\qquad (n\rightarrow\infty),
\end{equation}
where the induction is with respect to the lexicographic order of $(r,s)$.

Note that the base case, namely $(0,s)$ with $s\geq 0$, is implied by (\ref{exp-mom-Tn}).

Next, by the above lemma ${\mathbb E}(X_{n}^{\underline{r}}T_n^s)$ satisfies a recurrence of the form (\ref{gen-rec}) with all terms in $\psi_n$ being of a smaller lexicographic order. Thus, by using the induction hypothesis,
\[
\psi_n\sim\frac{4r}{n^2}{\mathbb E}(X_n^{\underline{r-1}} T_n^{s+1})+\lambda_{s-1,0}{\mathbb E}(X_n^{\underline{r}} T_n^{s-1})\sim\frac{8r+7s}{14^r7^s}n^{s-1}.
\]
The induction claim follows now from this by applying Lemma~\ref{asymp-transfer}.

Finally, the claimed Poisson limit law for $X_n$ follows from (\ref{main-claim}).\qed

\begin{figure}[!t]
    \centering
    \includegraphics[scale=0.9]{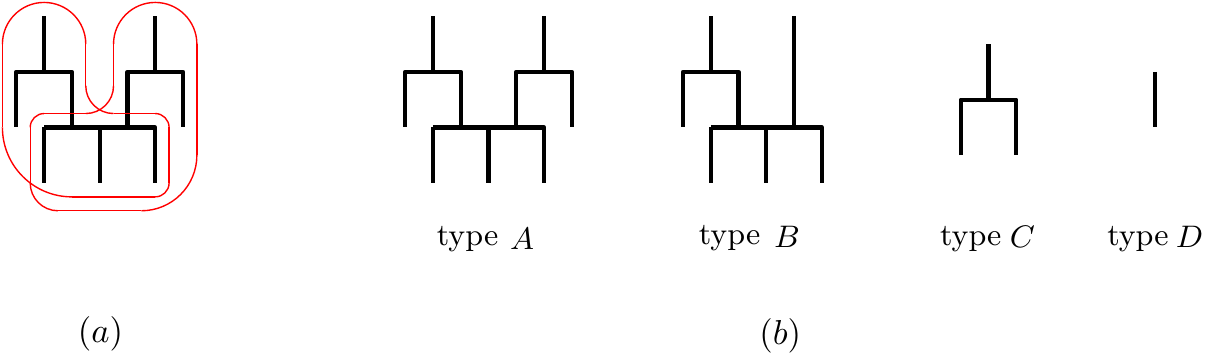}
    \caption{(a) The pattern of height $3$ which contains two overlapping patterns from Figure~\ref{height-2}-(b-i); (b) The types of patterns considered in the proof of the Poisson limit law for the pattern in Figure~\ref{height-2}-(b-i); in order that every external lineage belongs exactly to one type, a pattern of type $B$ is not allowed to be contained in a pattern of type $A$; also, the lineage in type $D$ is not an external lineage in $A$, $B$ or $C$.}\label{pattern-(b-i)}
\end{figure}

Next, we turn to the remaining pattern from Figure~\ref{height-2}-(b), namely the pattern in Figure~\ref{height-2}-(b-i). This pattern is different from the other four that satisfy a Poisson limit law because two occurrences of this pattern in a ranked tree-child network might overlap; see Figure~\ref{pattern-(b-i)}-(a). (For the pattern in this figure, we do not distinguish the ranks of the first two events.) Thus, we now need to consider patterns of types $A$-$D$ in order to set up the Markov chain for proving the Poisson limit law for this pattern; see Figure~\ref{pattern-(b-i)}-(b).

Assume that the number of patterns of type $A$, type $B$ and type $C$ in a ranked tree-child network with $n$ leaves is given by $a, b$ and $c$, respectively. Then, the changes and probabilities when adding one more event in the forward construction are listed in Table~\ref{bra-(b-i)} and Table~\ref{ret-(b-i)}.

\vspace{0.35cm}
\begin{table}[!h]
\begin{center}
\begin{tabular}{c|c|c|c|c}
   & type $A$ & type $B$ & type $C$ & probability \\ \hline\hline
   \multirow{2}{*}{$A$} & $-1$ & $0$ & $+1$ & $3a/n^2$ \\
   & $-1$ & $+1$ & $+1$ &$2a/n^2$ \\ \hline
   $B$ & $0$ & $-1$ & $+1$ & $4b/n^2$ \\ \hline
   $C$ & $0$ & $0$ & $0$ & $2c/n^2$ \\ \hline
   $D$ & $0$ & $0$ & $+1$ & $(n-5a-4b-2c)/n^2$
\end{tabular}
\end{center}
\vspace*{-.2cm}
\caption{The change of the number of patterns of type $A$, type $B$ and type $C$ (see Figure~\ref{pattern-(b-i)}-(b)) when the next event is a branching event.}\label{bra-(b-i)}
\end{table}

Using the transition probabilities from Table~\ref{bra-(b-i)} and Table~\ref{ret-(b-i)}, we obtain now the following lemma.

\begin{lmm}\label{rec-mom-(c-i)}
Denote by $Y_n$ and $\tilde{X}_n$ the number of occurrences of patterns of type $A$ and type $B$, respectively, in a random ranked tree-child network with $n$ leaves. Then, for $r,s,t\geq 0$, we have
\begin{align}
{\mathbb E}(Y_{n+1}^{\underline{r}}\tilde{X}_{n+1}^{\underline{s}}C_{n+1}^{\underline{t}})=\left(1-\frac{5r+4s+2t}{n}\right)^2&{\mathbb E}(Y_{n}^{\underline{r}}\tilde{X}_{n}^{\underline{s}}C_{n}^{\underline{t}})+\frac{t}{n}{\mathbb E}(Y_{n}^{\underline{r}}\tilde{X}_{n}^{\underline{s}}C_{n}^{\underline{t-1}})+\frac{4s}{n}{\mathbb E}(Y_{n}^{\underline{r+1}}\tilde{X}_{n}^{\underline{s-1}}C_{n}^{\underline{t}})\nonumber\\
&+\frac{4s}{n}{\mathbb E}(Y_{n}^{\underline{r}}\tilde{X}_{n}^{\underline{s-1}}C_{n}^{\underline{t+1}})+\frac{R_n}{n^2},\label{rec-(b-i)}
\end{align}
where $R_n$ is given by
\begin{align*}
t&(2-5r-4s-2t){\mathbb E}(Y_{n}^{\underline{r}}\tilde{X}_{n}^{\underline{s}}C_{n}^{\underline{t-1}})+4r{\mathbb E}(Y_{n}^{\underline{r-1}}\tilde{X}_{n}^{\underline{s}}C_{n}^{\underline{t+2}})-2s(1+10r+8s+4t){\mathbb E}(Y_{n}^{\underline{r+1}}\tilde{X}_{n}^{\underline{s-1}}C_{n}^{\underline{t}})\\
&+2st{\mathbb E}(Y_{n}^{\underline{r+1}}\tilde{X}_{n}^{\underline{s-1}}C_{n}^{\underline{t-1}})-8s{\mathbb E}(Y_{n}^{\underline{r}}\tilde{X}_{n}^{\underline{s-1}}C_{n}^{\underline{t+2}})+4s(2-5r-4s-2t){\mathbb E}(Y_{n}^{\underline{r}}\tilde{X}_{n}^{\underline{s-1}}C_{n}^{\underline{t+1}})\\
&+4s(s-1){\mathbb E}(Y_{n}^{\underline{r+2}}\tilde{X}_{n}^{\underline{s-2}}C_{n}^{\underline{t}})+8s(s-1){\mathbb E}(Y_{n}^{\underline{r+1}}\tilde{X}_{n}^{\underline{s-2}}C_{n}^{\underline{t+1}}).
\end{align*}
\end{lmm}

\begin{table}[!h]
\begin{center}
\begin{tabular}{c|c|c|c|c}
   & type $A$ & type $B$ & type $C$ & probability \\ \hline\hline
   $A$ & $-1$ & $0$ & $0$ & $20an^2$ \\ \hline
   \multirow{3}{*}{$A\ \& \ A$} & $-2$ & $0$ & $0$ & $9a(a-1)/n^2$ \\
   & $-2$ & $+1$ & $0$ & $12a(a-1)/n^2$ \\
   & $-2$ & $+2$ & $0$ & $4a(a-1)/n^2$ \\ \hline
   $B$ & $0$ & $-1$ & $0$ & $12b/n^2$ \\ \hline
   $B \&\ B$ & $0$ & $-2$ & $0$ & $16b(b-1)/n^2$ \\ \hline
   $C$ & $0$ & $0$ & $-1$ & $2c/n^2$ \\ \hline
   $C\ \&\ C$ & $+1$ & $0$ & $-2$ & $4c(c-1)/n^2$\\ \hline
   $D\ \&\ D$ & $0$ & $0$ & $0$ & $(n-5a-4b-2c)(n-5a-4b-2c-1)/n^2$ \\ \hline
   \multirow{2}{*}{$A\ \&\ B$} & $-1$ & $-1$ & $0$ & $24ab/n^2$\\
   & $-1$ & $0$ & $0$ & $16ab/n^2$ \\ \hline
   \multirow{2}{*}{$A\ \&\ C$} & $-1$ & $+1$ & $-1$ & $12ac/n^2$ \\
   & $-1$ & $+2$ & $-1$ & $8ac/n^2$ \\ \hline
   \multirow{2}{*}{$A\ \&\ D$} & $-1$ & $0$ & $0$ & $6a(n-5a-4b-2c)/n^2$ \\
   & $-1$ & $+1$ & $0$ & $4a(n-5a-4b-2c)/n^2$ \\ \hline
   $B\ \&\ C$ & $0$ & $0$ & $-1$ & $16bc/n^2$ \\ \hline
   $B\ \&\ D$ & $0$ & $-1$ & $0$ & $8b(n-5a-4b-2c)/n^2 $\\ \hline
   $C\ \&\ D$ & $0$ & $+1$ & $-1$ & $4c(n-5a-4b-2c)/n^2$
\end{tabular}
\end{center}
\vspace*{-.2cm}
\caption{The change of the number of patterns of type $A$ and type $B$ (see Figure~\ref{pattern-(b-i)}-(b)) when the next event is a reticulation event.}\label{ret-(b-i)}
\end{table}

\pf Using the results from Table~\ref{bra-(b-i)} and Table~\ref{ret-(b-i)}, ${\mathbb E}(Y_{n+1}^{\underline{r}}\tilde{X}_{n+1}^{\underline{s}}C_{n+1}^{\underline{t}})$ given $Y_n, \tilde{X}_n$ and $C_n$ can be written as a sum of terms of the form
\[
(Y_{n}+k)^{\underline{r}}(\tilde{X}_{n}+\ell)^{\underline{s}}(C_{n}+m)^{\underline{t}}
\]
which are multiplied with the probabilities. Here, $k,\ell$ and $m$ are suitable integers, e.g., for the contribution from the second sub-row of column $A\ \&\ C$ in Table~\ref{ret-(b-i)}, we have $k=-1,\ell=2,m=-1$. Then, we rewrite, e.g., this term as
\[
(Y_{n}-1)^{\underline{r}}(\tilde{X}_{n}+2)^{\underline{s}}(C_{n}-1)^{\underline{t}}=\frac{(Y_n-r)(\tilde{X}_n+2)(\tilde{X}_n+1)(C_n-t)}{Y_n(\tilde{X}_n-s)(\tilde{X}_n-s-1)C_n}Y_{n}^{\underline{r}}\tilde{X}_{n}^{\underline{s}}C_{n}^{\underline{t}}
\]
and similar for the other terms. The rest of the proof is just a long computation (which is best done with the help of Maple).\qed

The last lemma implies the following result which contains the Poisson limit law for the pattern in Figure~\ref{height-2}-(b-i).

\begin{pro}\label{limit-laws-(b-i)}
\begin{itemize}
\item[(a)] Let $Y_n$ be the number of occurrences of the pattern from Figure~\ref{pattern-(b-i)}-(a) in a random ranked tree-child network with $n$ leaves. Then,
\[
Y_n\stackrel{L_1}{\longrightarrow} 0,\qquad(n\rightarrow\infty).
\]
\item[(b)] Let $X_n$ be the number of occurrences of the pattern from Figure~\ref{height-2}-(b-i) in a random ranked tree-child network with $n$ leaves. Then,
\[
(X_n,C_n)\stackrel{d}{\longrightarrow}(X,C),\qquad (n\rightarrow\infty),
\]
where $X$ and $C$ are independent Poisson random variables with parameters $1/8$ and $1/4$, respectively.
\end{itemize}
\end{pro}

\pf We use induction with respect to the lexicographic order on $(s,r,t)$ to show that for all $r,s,t\geq 0$:
\begin{equation}\label{ind-claim}
{\mathbb E}(Y_n^{\underline{r}}\tilde{X}_n^{\underline{s}}C_n^{\underline{t}})\sim0^{r}\cdot\left(\frac{1}{8}\right)^s\cdot\left(\frac{1}{4}\right)^t,
\end{equation}
where $\tilde{X}_n$ is the random variable from Lemma~\ref{rec-mom-(c-i)} and we use the convention that $0^0:=1$. Note that the second, third and fourth term on the right hand side of (\ref{rec-(b-i)}) and all terms in $R_n$ are of a smaller lexicographic order than $(s,r,t)$. Also note that (\ref{rec-(b-i)}) has the form (\ref{gen-rec}).

Now, first the induction base holds because of (\ref{exp-mom-Cn}). Next, assume that the claim holds for all sequences $(s',r',t')$ which are lexicographic smaller than $(s,r,t)$. In order to prove the claim for $(s,r,t)$, we make a case distinction.

First, if $r>0$, then the induction hypothesis implies that (\ref{rec-(b-i)}) satisfies (\ref{gen-rec}) with $\psi_n=o(1/n)$. Thus, from Lemma~\ref{asymp-transfer}, we obtain that ${\mathbb E}(Y_n^{\underline{r}}\tilde{X}_n^{\underline{s}}C_n^{\underline{t}})=o(1)$ which proves the claim in this case.

Secondly, if $r=0$, then again by the induction claim, (\ref{rec-(b-i)}) is of the form (\ref{gen-rec}) with
\[
\psi_n\sim\frac{t}{n}{\mathbb E}(\tilde{X}_n^{\underline{s}}C_n^{\underline{t-1}})+\frac{4s}{n}{\mathbb E}(\tilde{X}_n^{\underline{s-1}}C_n^{\underline{t+1}})\sim\frac{8s+4t}{8^s4^tn}.
\]
Thus, by Lemma~\ref{asymp-transfer}
\[
{\mathbb E}(\tilde{X}_n^{\underline{s}}C_n^{\underline{t}})\sim\frac{1}{8^s4^t}.
\]
which also proves the claim in this case.

Finally, observe that (\ref{ind-claim}) implies part (a). Moreover, it also implies that
\[
(\tilde{X}_n,C_n)\stackrel{d}{\longrightarrow}(\tilde{X},C),
\]
where $\tilde{X}$ and $C$ are independent Poisson random variables with parameters $1/8$ and $1/4$, respectively. Part (b) follows from this by part (a), the relation $X_n=\tilde{X}_n+2Y_n$ and Slutsky's theorem.\qed

\paragraph{Normal Limit Laws.} In this paragraph, we consider the patterns in Figure~\ref{height-2}-(c) which according to Theorem~\ref{main-result-2} are both normal distributed. We will only give details for the pattern in Figure~\ref{height-2}-(c-i), the other pattern in Figure~\ref{height-2}-(c) is treated similarly.

\begin{figure}[!t]
    \centering
    \includegraphics[scale=0.9]{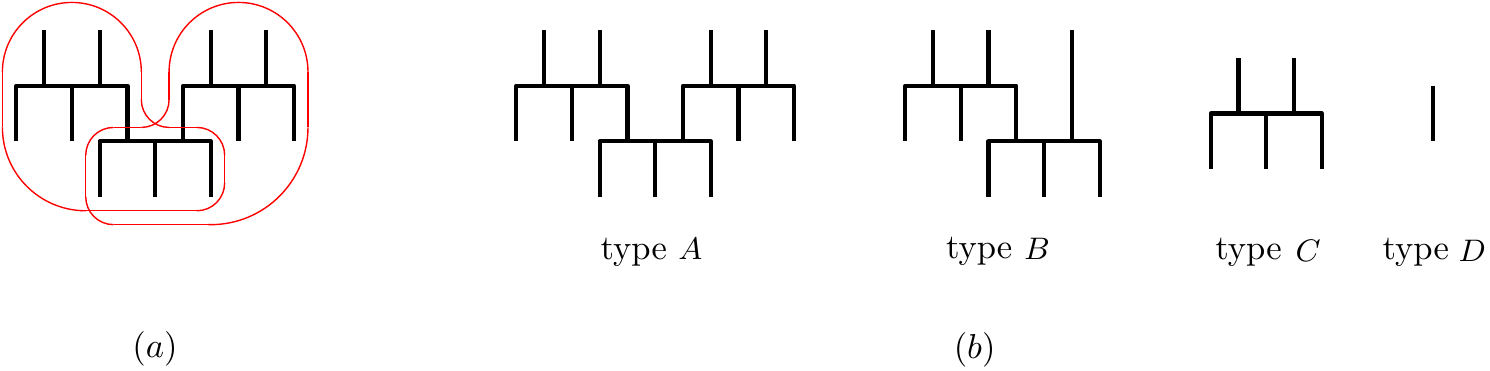}
    \caption{(a) The pattern of height $3$ which contains two overlapping patterns from Figure~\ref{height-2}-(c-i); (b) The types of pattern considered in the proof of the normal limit law for the pattern in Figure~\ref{height-2}-(c-i); in order that every external lineage belongs exactly to one type, a pattern of type $B$ resp. $C$ is not allowed to be contained in a pattern of type $A$ resp. $B$ and $A$; also, the lineage in type $D$ is not an external lineage in $A$, $B$ or $C$.}\label{pattern-(c-i)}
\end{figure}

First note that two of the patterns in Figure~\ref{height-2}-(c-i) again can overlap in a ranked tree-child network; see  Figure~\ref{pattern-(c-i)}-(a). (We do not distinguish the ranks of the first two events in this figure.) Thus, similar to the pattern in Figure~\ref{height-2}-(b-i), we will consider four types of patterns in the sequel; see Figure~\ref{pattern-(c-i)}-(b). Assume that a ranked tree-child network with $n$ leaves contains $a, b$ and $c$ patterns of type $A, B$ and $C$, respectively. We again have to list the changes and probabilities if one more event is attached in the forward construction, however, this time there are many cases. We list all of them in Tables~\ref{bra-(c-i)}-\ref{ret-(c-i)}.

Using these tables, we can first derive the means of the number of occurrences of the patterns from Figure~\ref{height-2}-(c-i) and Figure~\ref{pattern-(c-i)}-(a).

\begin{lmm}\label{mean-(c-i)}
Denote by $X_n$ and $Y_n$ the number of occurrences of the patterns from Figure~\ref{height-2}-(c-i) and Figure~\ref{pattern-(c-i)}-(a), respectively. Let $\rho_n:={\mathbb E}(X_n)$ and $\tau_n:={\mathbb E}(Y_n)$. Then,
\[
\rho_{n+1}=\left(1-\frac{5}{n}\right)^2\rho_n+\frac{4}{n}{\mathbb E}(T_n)-\frac{12}{n^2}{\mathbb E}(T_n)
\]
and
\[
\tau_{n+1}=\left(1-\frac{7}{n}\right)^2\tau_n-\frac{4}{n^2}{\mathbb E}(T_n)+\frac{4}{n^2}{\mathbb E}(T_n^2).
\]
The solutions of the above two recurrences are given by
\begin{equation}\label{exp-rho}
\rho_n=\frac{(1080n^5-16668n^4+96992n^3-261735n^2+319471n-135654)n}{20790(n-1)(n-2)(n-3)(n-4)(n-5)},\qquad (n\geq 6)
\end{equation}
and, for $n\geq 7$,
\begin{equation}\label{exp-tau}
\scalebox{1.2}{$\tau_n=\frac{2(4290n^7-125730n^6+1509970n^5-9550275n^4+33968326n^3-66128140n-24510098)n}{1576575(n-1)(n-2)(n-3)(n-4)(n-5)(n-6)(n-7)}$.}
\end{equation}
\end{lmm}

\pf Assume that $Y_n=d, X_n=e$ and $T_n=f$. Then, we have $d$ patterns of type $A$, $e-2d$ patterns of type $B$ and $f+e-d$ patterns of type $C$.

The recurrence for $\mu_n$ and $\rho_n$ are derived by considering the transition probabilities from the Tables~\ref{bra-(c-i)}-\ref{ret-(c-i)}. Note that both recurrences are of type (\ref{gen-rec}) whose exact solution is given by (\ref{sol-gen-rec}). Using the result from Lemma~\ref{mean-BiLaSt}, a corresponding result for ${\mathbb E}(T_n^2)$ (best computed with Maple) and straightforward computations (again best done with Maple) gives the claimed results for $\rho_n$ and $\tau_n$.\qed

\vspace*{0.2cm}
\begin{table}[!h]
\begin{center}
\begin{tabular}{c|c|c|c|c}
   & type $A$ & type $B$ & type $C$ & probability \\ \hline\hline
   \multirow{2}{*}{$A$} & $-1$ & $0$ & $0$ & $3a/n^2$ \\
   & $-1$ & $+1$ & $0$ & $4a/n^2$\\ \hline
   \multirow{2}{*}{$B$} & $0$ & $-1$ & $0$ & $3b/n^2$ \\
   & $0$ & $-1$ & $+1$ & $2b/n^2$ \\ \hline
   $C$ & $0$ & $0$ & $-1$ & $3c/n^2$ \\ \hline
   $D$ & $0$ & $0$ & $0$ & $(n-7a-5b-3c)/n^2$
\end{tabular}
\end{center}
\vspace*{-.2cm}
\caption{The change of the number of patterns of type $A, B$ and $C$ (see Figure~\ref{pattern-(c-i)}-(b)) when the next event is a branching event.}\label{bra-(c-i)}
\end{table}

\begin{table}[!t]
\begin{center}
\begin{tabular}{c|c|c|c|c}
   & type $A$ & type $B$ & type $C$ & probability \\ \hline\hline
   \multirow{4}{*}{$A$} & $-1$ & $0$ & $+1$ & $14a/n^2$\\
   & $-1$ & $0$ & $+2$ & $8a/n^2$ \\
   & $-1$ & $+1$ & $0$ & $16a/n^2$ \\
   & $-1$ & $+1$ & $+1$ & $4a/n^2$ \\ \hline
   \multirow{6}{*}{$A\ \&\ A$} & $-1$ & $0$ & $0$ & $4a(a-1)/n^2$ \\
   & $-2$ & $0$ & $+1$ & $a(a-1)/n^2$ \\
   & $-2$ & $+1$ & $0$ & $4a(a-1)/n^2$ \\
   & $-2$ & $+1$ & $+1$ & $8a(a-1)/n^2$ \\
   & $-2$ & $+2$ & $0$ & $16a(a-1)/n^2$ \\
   & $-2$ & $+2$ & $+1$ & $16a(a-1)/n^2$ \\ \hline
   \multirow{3}{*}{$B$} & $0$ & $0$ & $0$ & $8b/n^2$\\
   & $0$ & $-1$ & $+1$ & $10b/n^2$ \\
   & $0$ & $-1$ & $+2$ & $2b/n^2$ \\ \hline
   \multirow{6}{*}{$B\ \&\ B$} & $0$ & $-1$ & $0$ & $4b(b-1)/n^2$ \\
   & $0$ & $-1$ & $+1$ & $8b(b-1)/n^2$ \\
   & $0$ & $-2$ & $+1$ & $b(b-1)/n^2$ \\
   & $0$ & $-2$ & $+2$ & $4b(b-1)/n^2$ \\
   & $0$ & $-2$ & $+3$ & $4b(b-1)/n^2$ \\
   & $+1$ & $-2$ & $0$ & $4b(b-1)/n^2$ \\ \hline
   $C$ & $0$ & $0$ & $0$ & $6c/n^2$ \\ \hline
   \multirow{3}{*}{$C\ \&\ C$} & $0$ & $0$ & $-1$ & $c(c-1)/n^2$ \\
   & $0$ & $+1$ & $-2$ & $4c(c-1)/n^2$ \\
   & $+1$ & $0$ & $-2$ & $4c(c-1)/n^2$ \\ \hline
   $D\ \&\ D$ & $0$ & $0$ & $+1$ & $(n-7a-5b-3c)(n-7a-5b-3c-1)/n^2$
\end{tabular}
\end{center}
\caption{The change of the number of patterns of type $A, B$ and $C$ (see Figure~\ref{pattern-(c-i)}-(b)) when the next event is a reticulation event which is attached to one or two patterns of type $X$ with $X\in\{A,B,C,D\}$.}
\end{table}

\begin{table}[!t]
\begin{center}
\begin{tabular}{c|c|c|c|c}
   & type $A$ & type $B$ & type $C$ & probability \\ \hline\hline
    \multirow{7}{*}{$A\ \&\ B$} & $0$ & $-1$ & $0$ & $8ab/n^2$ \\
    & $-1$ & $-1$ & $+1$ & $2ab/n^2$ \\
    & $-1$ & $-1$ & $+2$ & $4ab/n^2$ \\
    & $-1$ & $0$ & $0$ & $8ab/n^2$ \\
    & $-1$ & $0$ & $+1$ & $16ab/n^2$ \\
    & $-1$ & $0$ & $+2$ & $16ab/n^2$ \\
    & $-1$ & $+1$ & $0$ & $16ab/n^2$ \\ \hline
    \multirow{5}{*}{$A\ \&\ C$} & $0$ & $0$ & $-1$ & $8ac/n^2$ \\
    & $-1$ & $0$ & $0$ & $2ac/n^2$ \\
    & $-1$ & $+1$ & $-1$ & $8ac/n^2$ \\
    & $-1$ & $+1$ & $0$ & $8ac/n^2$ \\
    & $-1$ & $+2$ & $-1$ & $16ac/n^2$ \\ \hline
    \multirow{3}{*}{$A\ \&\ D$} & $-1$ & $0$ & $+1$ & $2a(n-7a-5b-3c)/n^2$ \\
    & $-1$ & $+1$ & $0$ & $4a(n-7a-5b-3c)/n^2$ \\
    & $-1$ & $+1$ & $+1$ & $8a(n-7a-5b-3c)/n^2$ \\ \hline
    \multirow{5}{*}{$B\ \&\ C$} & $0$ & $-1$ & $0$ & $2bc/n^2$ \\
    & $0$ & $-1$ & $+1$ & $4bc/n^2$ \\
    & $0$ & $0$ & $-1$ & $8bc/n^2$ \\
    & $0$ & $0$ & $0$ & $8bc/n^2$ \\
    & $+1$ & $-1$ & $-1$ & $8bc/n^2$ \\ \hline
    \multirow{3}{*}{$B\ \&\ D$} & $0$ & $-1$ & $+1$ & $2b(n-7a-5b-3c)/n^2$ \\
    & $0$ & $-1$ & $+2$ & $4b(n-7a-5b-3c)/n^2$ \\
    & $0$ & $0$ & $0$ & $4b(n-7a-5b-3c)/n^2$ \\ \hline
    \multirow{2}{*}{$C\ \&\ D$} & $0$ & $0$ & $0$ & $2c(n-7a-5b-3c)/n^2$ \\
    & $0$ & $+1$ & $-1$ & $4c(n-7a-5b-3c)/n^2$
\end{tabular}
\end{center}
\vspace*{-.2cm}
\caption{The change of the number of patterns of type $A, B$ and $C$ (see Figure~\ref{pattern-(c-i)}-(b)) when the next event is a reticulation event which is attached to a pattern of type $X$ and a pattern of type $Y$ with $(X,Y)\in\{(A,B),(A,C),(A,D),(B,C),(B,D),(C,D)\}$.}\label{ret-(c-i)}
\end{table}

We next shift the means and derive the recurrence for the mixed moments.

\begin{lmm}\label{rec-normal}
Denote by $X_n$ and $Y_n$ the number of occurrences of the patterns from Figure~\ref{height-2}-(c-i) and Figure~\ref{pattern-(c-i)}-(a), respectively, in a random ranked tree-child network with $n$ leaves. Moreover, set $\mu_n:={\mathbb E}(T_n), \rho_n:={\mathbb E}(X_n), \tau_n:={\mathbb E}(Y_n)$ and $\bar{T}_n:=T_n-\mu_n, \bar{X}_n:=X_n-\rho_n, \bar{Y}_n:=Y_n-\tau_n$. Then, for all $r,s,t\geq 0$, we have
\begin{equation}\label{rec-(c-i)}
{\mathbb E}(\bar{Y}_{n+1}^{r}\bar{X}_{n+1}^{s}\bar{T}_{n+1}^{t})=\left(1-\frac{7r+5s+3t}{n}\right)^{2}{\mathbb E}(\bar{Y}_n^{r}\bar{X}_n^{s}\bar{T}_n^{t})+R_n
\end{equation}
with
\begin{equation}\label{toll-sequences}
R_n=\sum_{(s',r',t')}{\mathbb E}(\bar{Y}_{n}^{r'}\bar{X}_{n}^{s'}\bar{T}_{n}^{t'})\Lambda_{r',s',t'}(n),
\end{equation}
where the sum runs over $(s',r',t')$ which are of a smaller lexicographic order than $(s,r,t)$ and $\Lambda_{r',s',t'}(n)$ admits the complete asymptotic expansion:
\begin{equation}\label{exp-Lambda}
\Lambda_{r',s',t'}(n)\sim\sum_{\ell=0}^{\infty}\frac{\lambda_{r',s',t',\ell}}{n^{\ell}},\qquad (n\rightarrow\infty).
\end{equation}
Moreover, all terms in (\ref{toll-sequences}) with $(r'+s'+t')/2-\ell\geq (r+s+t)/2-1$ are given by
\begin{align}
\frac{4s}{n}&{\mathbb E}(\bar{Y}_{n}^{r}\bar{X}_{n}^{s-1}\bar{T}_{n}^{t+1})+\frac{8r}{7n}{\mathbb E}(\bar{Y}_{n}^{r-1}\bar{X}_{n}^{s}\bar{T}_{n}^{t+1})\nonumber\\
&+\binom{r}{2}\frac{80092}{540225}{\mathbb E}(\bar{Y}_{n}^{r-2}\bar{X}_{n}^{s}\bar{T}_{n}^{t})+\binom{s}{2}\frac{21916}{29645}{\mathbb E}(\bar{Y}_{n}^{r}\bar{X}_{n}^{s-2}\bar{T}_{n}^{t})+\binom{t}{2}\frac{24}{49}{\mathbb E}(\bar{Y}_{n}^{r}\bar{X}_{n}^{s}\bar{T}_{n}^{t-2})\nonumber\\
&-\frac{128}{539}st{\mathbb E}(\bar{Y}_{n}^{r}\bar{X}_{n}^{s-1}\bar{T}_{n}^{t-1})-\frac{32}{343}rt{\mathbb E}(\bar{Y}_{n}^{r-1}\bar{X}_{n}^{s}\bar{T}_{n}^{t-1})+\frac{712}{3773}rs{\mathbb E}(\bar{Y}_{n}^{r-1}\bar{X}_{n}^{s-1}\bar{T}_{n}^{t}).\label{main-terms}
\end{align}
\end{lmm}

\pf Assume that $Y_n, X_n$ and $T_n$ are given. Then, using the cases listed in the Tables~\ref{bra-(c-i)}-\ref{ret-(c-i)}, we can write the conditional expectation of $\bar{Y}_{n+1}^{r}\bar{X}_{n+1}^{s}\bar{T}_{n+1}^{t}$ given $Y_n, X_n$ and $T_n$ as a sum of terms of the form
\begin{equation}\label{shape-terms}
(\bar{Y}_n+\tau_n-\tau_{n+1}+k)^{r}(\bar{X}_n+\rho_n-\rho_{n+1}+\ell)^{s}(\bar{T}_n+\mu_n-\mu_{n+1}+m)^{t}
\end{equation}
multiplied with the probabilities from the Tables~\ref{bra-(c-i)}-\ref{ret-(c-i)} (where $a=\bar{Y}_n+\tau_n, b=\bar{X}_n+\rho_n-2(\bar{Y}_n+\tau_n), c=\bar{T}_n+\mu_n+\bar{X}_n+\rho_n-\bar{Y}_n-\tau_n$). Here, $k,\ell,m$ are integers depending on the case considered. For convenience, we will call this sum $S$ in the sequel.

Now, use the expansion:
\begin{equation}\label{bin-th}
(\bar{X}_n+\rho_n-\rho_{n+1}+\ell)^{s}=\sum_{j=0}^{s}\binom{s}{j}\bar{X}_n^{j}(\rho_n-\rho_{n+1}+\ell)^{s-j}.
\end{equation}
First, by replacing the middle term in (\ref{shape-terms}) by $\bar{X}_n^{s}$ (the first term in the expansion (\ref{bin-th})) and using Maple to sum up all these terms multiplied with their probabilities (which contain at most $\bar{X}_n^2$), we find that $\bar{X}_n^{s+2}$ and $\bar{X}_n^{s+1}$ do not appear in this sum. Likewise, by replacing the middle term in  (\ref{shape-terms}) by $s\bar{X}_n^{s-1}(\rho_n-\rho_{n+1}+\ell)$ (the second term in the expansion (\ref{bin-th})), we see that $\bar{X}_n^{s+1}$ does not appear. Thus, the highest power of $\bar{X}_n$ in the sum $S$ is $\bar{X}_n^{s}$. Next, by collecting these highest terms and repeating the above argument with the first term in (\ref{shape-terms}), we see that terms with $\bar{Y}_n^{r+2}\bar{X}_n^{s}$ and $\bar{Y}_n^{r+1}\bar{X}_n^s$ do not occur in $S$. Finally, a similar line of reasoning shows that terms with $\bar{Y}_n^r\bar{X}_n^{s}\bar{T}_n^{t+2}$ and $\bar{Y}_n^r\bar{X}_n^{s}\bar{T}_n^{t+1}$ do not appear in $S$ as well. Overall, this shows that sum in (\ref{toll-sequences}) is over the indicated range.

Next, the claimed expansion for $\Lambda_{r',s',t'}(n)$ follows by expanding $\tau_{n}-\tau_{n+1}, \rho_n-\rho_{n+1}$ and $\mu_n-\mu_{n+1}$ (see (\ref{exp-mu}), (\ref{exp-rho}) and (\ref{exp-tau}), respectively; note that these expansions are all of the form (\ref{exp-Lambda})) and pointing out that the probabilities might contain factors of the form $\tau_n^2, \rho_n^2, \mu_n^2, \tau_n\rho_n, \tau_n\mu_n,$ or $\rho_n\mu_n$ which are however divided by $n^2$ and thus expanding them gives also terms of the form (\ref{exp-Lambda}).

Finally, in order to find all terms with $(r'+s'+t')/2-\ell\geq (r+s+t)/2-1$, we proceed as follows: if we expand the first factor in (\ref{shape-terms}) and retain only the terms which contain $\bar{Y}_n^r, \bar{Y}_n^{r-1}, \bar{Y}_n^{r-2}, \bar{Y}_n^{r-3}$ and ${\bar Y}_n^{r-4}$, then we see that we only loose terms $\bar{Y}_n^{r'}\bar{X}_n^{s'}\bar{T}_n^{t'}$ in $S$ with
\[
(r'+s'+t')/2-\ell\leq (r'+s'+t')/2\leq (r-5+s+t+2)/2<(r+s+t)/2-1.
\]

Similarly, we only loose smaller order terms when we just keep the terms with $\bar{X}_n^{s},\ldots,\bar{X}_n^{s-4}$ and the terms with $\bar{T}_n^t,\ldots,\bar{T}_n^{t-4}$ in the expansion of the second and third factor in (\ref{shape-terms}), respectively. Thus, we only need to retain a fixed number (which does not depend on $r,s$ and $t$) of terms in each of the terms of (\ref{shape-terms}) from $S$. Then, the rest of the computation can then be done with Maple.\qed

The recurrence from the above lemma can now be used to prove the normal limit law of the pattern in Figure~\ref{height-2}-(c-i) (which completes the proof of Theorem~\ref{main-result-2}). In fact, we have a more general result.

\begin{pro}\label{clt-(c-i)}
Let the notation be as in Lemma~\ref{rec-normal}. Then, as $n\rightarrow\infty$,
\[
\frac{1}{n}(Y_n-{\mathbb E}(Y_n), X_n-{\mathbb E}(X_n), T_n-{\mathbb E}(T_n))\stackrel{d}{\longrightarrow}N({\bf 0},\Sigma),
\]
where $N({\bf 0},\Sigma)$ denotes a trivariate normal distribution with mean vector ${\bf 0}$ and covariance matrix
\[
\Sigma=\begin{pmatrix}\frac{1002796}{203664825} & \frac{433528}{62537475} & -\frac{32}{13377} \\[6pt] \frac{433528}{62537475} & \frac{4575916}{137582445} & -\frac{608}{119119} \\[6pt] -\frac{32}{13377} & -\frac{608}{119119} & \frac{24}{637}
\end{pmatrix}.
\]
\end{pro}

\vspace*{0.35cm}\noindent{\it Proof of Proposition~\ref{clt-(c-i)}.} We use induction with respect to the lexicographic order of $(s,r,t)$ to show that for all $r,s,t\geq 0$:
\[
{\mathbb E}(\bar{Y}_n^{r}\bar{X}_n^{s}\bar{T}_n^t)\sim c_{r,s,t}n^{(r+s+t)/2},
\]
where $c_{r,s,t}:={\mathbb E}(N_1^rN_2^sN_3^t)$ with $N=(N_1,N_2,N_3)$ the trivariate normal distribution with mean vector ${\bf 0}$ and covariance matrix $\Sigma$.

First, the base case is implied by Proposition~\ref{mom-Tn}. Thus, we may assume that the claim holds for all sequences $(s',r',t')$ which are lexicographic smaller than $(s,r,t)$.

In order to prove the claim for $(s,r,t)$, we observe that (\ref{rec-(c-i)}) has the form (\ref{gen-rec}) with the sequence $\psi_n$ satisfying $\psi_n\sim \tilde{c}_{r,s,t}n^{(r+s+t)/2-1}$ where
\begin{align}
\tilde{c}_{r,s,t}:=&4sc_{r,s-1,t+1}+\frac{8r}{7}c_{r-1,s,t+1}+\binom{r}{2}\frac{80092}{540225}c_{r-2,s,t}+\binom{s}{2}\frac{21916}{29645}c_{r,s-2,t}+\binom{t}{2}\frac{24}{49}c_{r,s,t-2}\nonumber\\ &\quad -\frac{128}{539}stc_{r,s-1,t-1}-\frac{32}{343}rtc_{r-1,s,t-1}+\frac{712}{3773}rsc_{r-1,s-1,t}.\label{tilde-c}
\end{align}
Now, recall that by Isserlis' theorem, we have
\[
c_{r,s-1,t+1}=\sum_{{\mathcal P}}\prod_{\{P_i,P_j\}\in {\mathcal P}}{\rm Cov}(P_i,P_j),
\]
where the sum runs over all partitions of the multiset $\{N_1,\ldots,N_1,N_2,\ldots,N_2,N_3,\ldots,N_3\}$ into pairs where $N_1$ is repeated $r$ times, $N_2$ is repeated $s-1$ times and $N_3$ is repeated $t+1$ times. By fixing one $N_3$ and considering its pairing with all other elements, we obtain that
\[
c_{r,s-1,t+1}=t\Sigma_{3,3}c_{r,s-1,t-1}+(s-1)\Sigma_{3,2}c_{r,s-2,t}+r\Sigma_{3,1}c_{r-1,s-1,t}.
\]
Moreover, by a similar argument for $c_{r-1,s,t+1}$, we have
\[
c_{r-1,s,t+1}=t\Sigma_{3,3}c_{r-1,s,t-1}+(s-1)\Sigma_{3,2}c_{r-1,s-1,t}+r\Sigma_{3,2}c_{r-2,s,t}.
\]
Plugging this into (\ref{tilde-c}), we obtain that
\begin{align}
\tilde{c}_{r,s,t}=&\binom{r}{2}\frac{1002796}{7022925}c_{r-2,s,t}+\binom{s}{2}\frac{4575916}{6551545}c_{r,s-2,t}+\binom{t}{2}\frac{24}{49}c_{r,s,t-2}\nonumber\\
&\quad-\frac{608}{7007}stc_{r,s-1,t-1}-\frac{32}{637}rtc_{r-1,s,t-1}+\frac{433528}{2501499}rsc_{r-1,s-1,t}.\label{ref-tilde-c}
\end{align}
Next, again by Isserlis' theorem, we have the recurrences:
\begin{align*}
c_{r,s,t}&=(r-1)\Sigma_{1,1}c_{r-2,s,t}+s\Sigma_{1,2}c_{r-1,s-1,t}+t\Sigma_{1,3}c_{r-1,s,t-1};\\
c_{r,s,t}&=r\Sigma_{2,1}c_{r-1,s-1,t}+(s-1)\Sigma_{2,2}c_{r,s-2,t}+t\Sigma_{2,3}c_{r,s-1,t-1};\\
c_{r,s,t}&=r\Sigma_{3,1}c_{r-1,s,t-1}+s\Sigma_{3,2}c_{r,s-1,t-1}+(t-1)\Sigma_{3,3}c_{r,s,t-2}.
\end{align*}
Multiplying the first by $29r/2$, the second by $21s/2$ and the third by $13t/2$ gives as right-hand side exactly (\ref{ref-tilde-c}). Thus,
\[
\tilde{c}_{r,s,t}=\left(\frac{29}{2}r+\frac{21}{2}s+\frac{13}{2}t\right)c_{r,s,t}.
\]
Overall, we have shown that (\ref{rec-(c-i)}) has the form (\ref{gen-rec}) with
\[
\psi_n\sim\left(\frac{29}{2}r+\frac{21}{2}s+\frac{13}{2}t\right)c_{r,s,t}n^{(r+s+t)/2-1}.
\]
The induction claim follows now from this by applying Lemma~\ref{asymp-transfer}. This completes the proof.\qed

\section{Conclusion}\label{con}

The main purpose of this paper was to study the number of occurrences of patterns on the fringe of ranked tree-child networks. More precisely, we strengthened the weak law of large numbers for the number of tridents from \cite{BiLaSt} by proving a central limit theorem (which completed the classification of the limit laws of patterns of height $1$) and studied the limit laws of all patterns of height $2$. There are only three possible limit laws: normal, Poisson and degenerate; accordingly, pattern either occur frequently (normal pattern), sporadically (Poisson pattern) or not all (degenerate pattern). We think that this behavior will persist for patterns of any height.

In fact, we can formulate a conjecture for the limit law of any fringe pattern which is defined as a connected substructure of a ranked tree-child network which has entirely evolved from a fixed set of lineages by consecutively adding branching and reticulation events. The limit law of such a pattern can be obtained recursively as follows.

\begin{conj}
Let $F$ be a fringe pattern. Denote by $P$ resp. $P_1$ and $P_2$ the patterns which are obtained from it by removing the last event. (Here, the second case is only possible if the last event is a reticulation event and the pattern gets disconnected when this event is removed.) Then, we have the following cases.
\begin{itemize}
\item[(a)] If $P$ is a normal pattern, then $F$ is a Poisson pattern; in all other cases for $P$, the pattern $F$ is a degenerate pattern.
\item[(b)] If $P_1$ and $P_2$ are both normal patterns, then $F$ is also a normal pattern; if $P_1$ is a normal pattern and $P_2$ is a Poisson pattern or vice versa, then $F$ is a Poisson pattern; in all remaining cases for $P_1$ and $P_2$, $F$ is a degenerate pattern.
\end{itemize}
\end{conj}

This conjecture is easily seen to be consistent with the limit laws for the patterns from Figure~\ref{height-2}; see Theorem~\ref{main-result-2}. It is also consistent with the limit law results for the two patterns of height $3$ above: the pattern from Figure~\ref{pattern-(b-i)}-(a) was shown to be degenerate (see Proposition~\ref{limit-laws-(b-i)}, (a)) and indeed the pattern splits into two Poisson patterns when the last event is removed; likewise, the pattern from Figure~\ref{pattern-(c-i)}-(a) is a normal pattern (see Proposition~\ref{clt-(c-i)}) and it splits into two normal patterns when the last event is removed.

For the proof of the results in this paper (Theorem~\ref{main-result-1} and Theorem~\ref{main-result-2}) we extended the approach from \cite{BiLaSt} which was based on coupling with a Markov chain and the method of moments. The method can be applied to prove further cases of the above conjecture; however, technical details become more and more demanding for patterns of increasing height. Thus, in order to proof the conjecture in its full generality, we think that a new (less computation-intensive) approach has to be devised.

We conclude by pointing out that one of the main contributions of the current study is that it constitutes the first such study for a class of phylogenetic networks. Indeed, it would be interesting to undertake similar studies for other classes of phylogenetic networks, e.g., the class of tree-child networks. However, for this class, almost nothing is known so far about the occurrence of patterns. In fact, the only result we are aware of was proved in \cite{DiSeWe} where the authors showed that the number of cherries is $o(n)$ for almost all tree-child networks with $n$ leaves. Thus, a (wild) guess would be that the limit law of this pattern is again the Poisson law; however, tools for proving this (as well as deriving limit laws for other patterns) are completely lacking. We leave this as an open problem.

\paragraph*{Acknowledgments.} We thank Yu-Sheng Chang and Guan-Ru Yu for joining the discussion on the research presented in this paper.

\nopagebreak

\end{document}